\magnification=1200
\input amstex
\documentstyle{amsppt}
\hoffset=-0.5pc
\nologo
\vsize=57.2truepc
\hsize=38.5truepc
\spaceskip=.5em plus.25em minus.20em

\define\Ho{\roman H}
\define\Bobb{\Bbb}
\define\etad{\eta}
\define\fra{\frak}
\define\smcc{}

\define\adamhilt{1}
\define\avrahalp{2}
\define\baueslem{3}
\define\berikash{4}
   \define\brown{5}
  \define\cartan{6}
\define\carteile{7}
    \define\chen{8}
\define\degrmosu{9}
 \define\doldone{10}
 \define\doldtwo{11}
\define\eilenber{12}
\define\eilenmac{13}
\define\gugenhei{14}
\define\gugentwo{15}
\define\gugenlam{16}
 \define\gulasta{17}
\define\gulstatw{18}
\define\gugenmay{19}
\define\gugenmil{20}
\define\gugenmun{21}
\define\gugensta{22}
\define\halperin{23}
\define\halpstas{24}
  \define\heller{25}
  \define\hilton{26}
  \define\habili{27}
\define\homotype{28}
\define\perturba{29}
\define\cohomolo{30}
\define\modpcoho{31}
\define\intecoho{32}
 \define\abelian{33}
\define\berikas{34}
\define\holomorp{35}
\define\omni{36}
\define\huebkade{37}
\define\huebstas{38}
\define\husmosta{39}
\define\kadeishv{40}
\define\kadeifou{41}
\define\kadeifiv{42}
\define\lambstas{43}
\define\liulevon{44}
\define\liulevtw{45}
 \define\maclane{46}
   \define\meyer{47}
\define\mooreone{48}
\define\mooretwo{49}
\define\moorethr{50}
\define\moorsmit{51}
\define\munkholm{52}
 \define\quillen{53}
\define\sanebone{54}
\define\sanebthr{55}
\define\sanebsev{56}
   \define\serre{57}
    \define\shih{58}
 \define\smirnov{59}
\define\stasheff{60}
\define\stashalp{61}
\define\sullivan{62}
\define\sullitwo{63}
    \define\wall{64}

\topmatter
\title
Minimal free multi models for chain algebras
\endtitle
\author Johannes Huebschmann
\endauthor
\affil
Universit\'e des Sciences et Technologies de Lille
\\
UFR de Math\'ematiques
\\
CNRS-UMR 8524
\\
F-59 655 VILLENEUVE D'ASCQ C\'edex, France
\\
Johannes.Huebschmann\@math.univ-lille1.fr
\endaffil
\date{May 10, 2004}
\enddate
\dedicatory To the memory of G. Chogoshvili
\enddedicatory
\abstract{ Let $R$ be a local ring and $A$ a connected
differential graded algebra over $R$ which is free as a graded
$R$-module. Using homological perturbation theory techniques, we
construct a minimal free multi model for $A$ having properties
similar to that of an ordinary minimal model over a field; in
particular the model is unique up to isomorphism of multialgebras.
The attribute \lq multi\rq\ refers to the category of
multicomplexes.}
\endabstract

\endkeywords
\endtopmatter
\document

\leftheadtext{Johannes Huebschmann}
\rightheadtext{Minimal free multi models for chain algebras}

\beginsection Introduction

Let $R$ be a commutative ring with 1, and let $A$ be a connected
differential graded algebra over $R$ which is free as a graded $R$-module,
endowed with the obvious augmentation map.
For example, $A$ could be the chains on the loop space $\Omega X$ of a
simply connected space $X$. As in
\cite\huebkade, we refer to a
differential graded algebra of the kind $(\roman T[V],d)$, where
$\roman T[V]$ denotes the graded tensor algebra on a free graded $R$-module $V$,
together with a morphism $(\roman T[V],d)\to(A,d)$ of differential graded
algebras which is also a chain equivalence, as a {\it free model for\/} $A$.
The approach in \cite\huebkade\ provides a small free model, and we recall
briefly the construction: Let $\Omega \roman B A$ be  the cobar construction
on the bar construction $\roman B A$, let $F{\Ho}(J\roman BA)$ be a free
resolution (in the category of $R$-modules) of the homology ${\Ho}(J\roman BA)$
of the coaugmentation coideal $J\roman BA$ of the bar construction $BA$, and
consider the tensor algebra $\roman T[s^{-1}F {\Ho}(J\roman BA)]$ on the
desuspension $s^{-1}F {\Ho}(J\roman BA)$ of $F {\Ho}(J\roman BA)$. A suitable
homological perturbation theory argument, applied to these data, enabled us to
construct
a differential $d$ on $\roman T[s^{-1}F {\Ho}(J\roman BA)]$ and a morphism
$$
(\roman T[s^{-1}F {\Ho}(J\roman BA)],d)
\to (\roman T[s^{-1}(J\roman BA)],d_{\Omega})
= \Omega \roman B A
$$
of differential graded algebras which is also chain equivalence; the composite
of this chain equivalence with the standard adjoint chain equivalence
$\Omega \roman B A \to A$ then yields a small free model for $A$.

In particular, when $R$ is a local ring which is as well a principal ideal
domain and when $F {\Ho}(J\roman BA)$ is a minimal resolution of the homology
${\Ho}(J\roman BA)$, the differential graded algebra
$(\roman T[s^{-1}F {\Ho}(J\roman BA)],d) $ together with the comparison map
into $A$ is what has been called a {\it minimal free model for\/} $A$
in \cite\huebkade. According to \cite\huebkade\ (5.11), such
a minimal free model exists and is unique up to isomorphism of chain algebras.
When the local ring is no longer a principal ideal domain, this approach still
yields a small free model but not a free minimal one in the naive sense,
cf. \cite\huebkade\ (5.12). In the present paper we shall
show that the resulting small free model is minimal as an algebra
in the category of {\it multicomplexes\/} or, equivalently, as a
{\it multialgebra\/} (precise definitions will be given in the next section)
and, given an augmented connected differential graded algebra $A$ that is free
as a module over the local ring $R$, we shall in fact establish existence and
uniqueness of what we shall call a {\it minimal free multi model\/} for $A$.
See Theorem 3.10 below for details. The idea of using this additional structure
is related with the more familiar one of using a filtration as an additional
piece of structure, cf. e.~g. \cite\halpstas.
Indeed, a multicomplex structure is equivalent to that of a filtered chain
complex having the property that the associated (bi)graded object
is free over the ground ring.
Multicomplexes occur at various places in the literature;
historical comments will be given in the next section.
A special case of a multicomplex arises from an ordinary chain complex
with the degree filtration, cf. (1.9) above.

Here is an outline of the contents of the paper. In Section 1 we recall the
concept of a multicomplex and introduce that of a {\it multialgebra\/}.
A special case of a
multialgebra is an ordinary differential graded algebra with the degree
filtration. We also introduce appropriate notions of morphism and of homotopy.
In Section 2 we explore free multialgebras, and in Section 3 we study minimal
free multialgebras over a local ring. In particular we shall show that, over
an arbitrary local ring, a differential graded algebra that is free as a module
over the ground ring, viewed as a multialgebra in the sense explained above,
has a minimal free model in the category of multialgebras that is unique up to
isomorphism. Details will be given in Theorem 3.10. Some comments about the
significance of this result and about its relationship with the literature
will be given in Remark 3.11.

The ground ring will be denoted by $R$ throughout, and graded and bigraded
modules  will always be free over the ground ring $R$ unless they are
explicitly specified otherwise; the notions of chain equivalence and weak
equivalence (i.~e. isomorphism on homology) are then equivalent, and we shall
use the term \lq weak equivalence\rq\ only when there is a difference between
the two. The same kind of remark applies to the concepts of
{\it multiequivalence\/} and {\it weak multiequivalence\/} introduced in
(1.11) and (1.12.1) below. The reader will have no trouble to replace
\lq free over $R$\rq\ with \lq projective over $R$\rq. We shall stick to the
free case to avoid unnecessary complications with language and terminology.
Our notation is the same as that in e.~g.
\cite\perturba, \cite\huebkade\ and \cite\munkholm.
Graded and bigraded algebras will always be assumed to be augmented.

This paper is dedicated to the memory of G. Chogoshvili. Within
the tradition on algebraic and topological research in Georgia
which goes back to him, the ideas which led to multicomplexes and
multialgebras are well represented, cf. e.~g.
\cite\berikash,\,\cite\huebkade,\,\cite\kadeishv--\cite\kadeifiv,\,
\cite\sanebone--\cite\sanebsev. This list is certainly not
exhaustive.

\beginsection 1. Multicomplexes and multialgebras

Let $R$ be a commutative ring with 1, taken henceforth as ground ring.
A {\it multicomplex\/} is a bigraded $R$-module together with a differential
on the associated graded module that preserves {\it column filtration\/}
(see Definition 1.4.1 below for details).
Taking components we arrive at the following.

\noindent
{\smc Definition 1.1.}
A {\it multicomplex\/} $X$ is a bigraded $R$-module
$\{X_{p,q}\}_{p,q \in \Bobb Z}$, together with $R$-linear morphisms
$$
d^i \colon X_{p,q}\longrightarrow X_{p-i,q+i-1},\quad i = 0,1,\cdots
$$
such that, for each $n \geq 0$,
$
\sum_{i+j=n} d^id^j = 0.
$

Henceforth we shall refer to $d = \{d^0, d^1, d^2, \dots \}$ as a
{\it multidifferential\/}. Notice that, for each $\ell$, the operator $d^0$
is a differential $X_{\ell,*} \to X_{\ell,*-1}$ but, for $j \geq 1$,
the operator $d^j$ is not necessarily a differential. We shall refer to
$d^0$ as the {\it vertical differential\/}. Likewise  we shall occasionally
refer to $d^1$ as a {\it horizontal operator\/}. When $d^0$ is zero,
for each $\ell$, the operator $d^1$ is manifestly a differential
$X_{*,\ell} \to X_{*-1,\ell}$:
we shall then refer to it as a {\it horizontal differential\/}.
A bicomplex may be viewed as a multicomplex with $d^i=0$ for $i \geq 2$.

The multicomplex terminology goes back at least to {\smcc
Liulevicius\/}~\cite\liulevon; without reference to an explicit
name, the structure has been exploited in  \cite\gugenmay,\
\cite\liulevtw,\ \cite\wall. A triangular complex in the sense of
\cite\heller\ is a special case of a multicomplex, and there is a
close relationship between multicomplexes and the predifferential
theory developed in \cite\berikash, cf. the proof of Theorem 3.10
below as well as \cite\sanebone,\,\cite\sanebthr. Multicomplexes
play a major role in {\it homological perturbation theory\/}, cf.
e.~g. Section 2 of \cite\modpcoho\ and Section 1 of
\cite\intecoho. More details and historical comments about
homological perturbation theory may be found e.~g. in
\cite\huebkade. A \lq\lq recursive structure of triangular
complexes\rq\rq, a concept isolated in Section 5 of \cite\heller,
is in fact an example of what was later identified as a
perturbation. In \cite\omni, certain algebraic structures behind
the spectral sequence of a foliation are explored by means of a
multialgebra version of the Maurer-Cartan algebra.

Given a bigraded $R$-module $X$, we shall refer to the graded $R$-module
$CX$, where
$$
CX_n = \sum_{i+j=n} X_{i,j},
$$
as the corresponding {\it total object\/}.

For a multicomplex $X$, the formal infinite sum $d = \sum d^j$
defines an operator on the total object $CX$ whenever the sum is finite in each
degree, and $(CX, d)$ is then a chain complex;
we refer to this situation by saying that
$(CX, d)$ is {\it well defined\/}.
This will manifestly be the case
when the column filtration is {\it bounded below\/} (cf.
e.g. \cite\maclane) in the sense that, for each degree $n$ (of $CX$),
there is an integer $s = s(n)$ such that
$$
X_{p,q} = 0\quad  \text{whenever} \quad p < s.
$$
Henceforth a multicomplex $X$ will be {\it assumed to be bounded
below\/} in this sense.

\noindent
{\smc Definition 1.2.}
Given a multicomplex $X$, its {\it total complex\/} is the chain complex
$(CX, d)$, where
$$
d = \sum d^j.
$$

A more rigorous description in the language of {\it assembly functors\/}
\cite\husmosta\ may be found in \cite\meyer.

An ordinary chain complex $C$ may be viewed as a multicomplex
in an obvious way. More precisely,
$$
C_{i,0} = C_i,\quad C_{i,\ne 0} =  0,
\quad d^1 = d,\quad d^j = 0 \text {\ for \ } j \ne 1,
\tag1.3
$$
yields a multicomplex whose total complex is just $C$.
We refer to (1.3) as the {\it associated multicomplex\/}.

\noindent
{\smc Definition 1.4.1.} The {\it column filtration\/} of a bigraded $R$-module
is the ascending filtration $\{F_p\}$ given by
$$
(F_p(X))_{i,j} = \cases
X_{i,j} &\text{\ if \ } i \leq p,\\
0 &\text{\ otherwise.}
\endcases
$$

\noindent
{\smc Definition 1.4.2.} The {\it row filtration\/} of a bigraded $R$-module
is the descending filtration $\{F^q\}$ given by
$$
(F^q(X))_{i,j} = \cases
X_{i,j} &\text{\ if \ } j \geq q,\\
0 &\text{\ otherwise.}
\endcases
$$

Given a multicomplex $X$, the row and column filtrations induce
corresponding filtrations on the total complex $CX$; in particular the
filtrations are compatible with the differential on the latter. We then refer
to these filtrations as {\it column\/} and {\it row filtrations\/} as well.
Moreover, the sum
$$
\partial = d^1 + d^2 + \dots
$$
is then what is called a {\it perturbation\/} of the differential $d^0$ on $CX$
with respect to the column filtration, that is, each $d^j$ lowers column
filtration by $j$.

\proclaim{Proposition 1.5}
Let $X$ be a bigraded $R$-module, and let $d$ be a differential on its total
object $CX$ that is compatible with the column filtration. Then the components
$$
d^i \colon X_{p,q} \longrightarrow X_{p-i,q+i-1},\quad i = 0,1,\cdots
$$
of $d$ endow $X$ with multicomplex structure in such a way that
totalization yields the original data.
\endproclaim

We shall need appropriate notions of morphism of multicomplexes and of
homotopy between such morphisms. To handle them concisely, we introduce the
following terminology; our description differs from the notions of morphism
given in \cite\liulevon,\ \cite\liulevtw, and \cite\meyer.

\noindent
{\smc Definition 1.6.}
Given two bigraded $R$-modules $X$ and $Y$, a {\it multimorphism of bigraded
$R$-modules of degree\/} $\etad$, written as $f \colon X \longrightarrow Y$,
consists of a sequence $f = \{f^k\}_{k \geq \ell}$ of $R$-module morphisms
$$
f^k \colon X_{p,q} \longrightarrow Y_{p-k,q+k+\etad}
$$
where $\ell$ is a (possibly negative) integer.
We refer to the $f^k$'s as the {\it components\/} of $f$,
and we denote the degree of $f$ by $|f|$ as usual.
We shall then write
$$
f= f^{\ell} + \dots + f^0 + f^1 + f^2 + \cdots \colon X \longrightarrow Y.
$$
Here the infinite sum is to be understood in a formal way. However, when $CX$
is well defined, this infinite sum converges in the sense that in each degree
only finitely many terms are non-zero. We note that a multimorphism
$f \colon X \longrightarrow Y$ preserves column filtrations
if and only if it is of the form
$$
f= f^0 + f^1 + f^2 + \cdots \colon X \longrightarrow Y.
$$

\noindent
{\smc Definition 1.7.1.}
Given two multimorphisms $f= f^{\ell} + \dots + f^0 + f^1 + f^2 +
\cdots \colon X \longrightarrow Y$ and
$g= g^{\ell'} + \dots + g^0 + g^1 + g^2 +
\cdots \colon Y \longrightarrow Z$ of bigraded modules
of degree $\etad$ and $\etad'$ respectively, the {\it composite\/}
$g \circ f$ is the  multimorphism $g \circ f \colon X \longrightarrow Z$
of bigraded $R$-modules of degree $\etad+\etad'$, where
$$
(g \circ f)^k = \sum_{i+j=k} g^if^j,
$$
that is, $g \circ f$ is obtained by a formal evaluation of the
\lq composition\rq
$$
(f^{\ell} + \dots + f^0 + f^1 + f^2 + \cdots)
(g^{\ell'} + \dots + g^0 + g^1 + g^2 + \cdots).
$$
This evaluation makes sense since, for each $k$, the sum $\sum_{i+j=k} g^if^j$
is finite. For example, given a bigraded $R$-module $X$, a multimorphism
$$
d = d^0 + d^1 + d^2 + \cdots \colon X \longrightarrow X
$$
of degree $-1$ yields a multicomplex structure on $X$ if and only if,
as a multimorphism of bigraded $R$-modules, the composite $d\circ d$ is zero.
This operation of composition of multimorphisms is plainly associative.
Henceforth we shall discard the symbol \lq $\circ$\rq \ and write
$g\,f = g \circ f$ etc.

\noindent
{\smc Remark 1.7.2.}
The bigraded $R$-modules together with a suitable choice of multimorphisms
constitute a category in an obvious fashion. In particular, invertible
multimorphisms of the kind
$f= f^0 + f^1 + f^2 + \cdots \colon X \longrightarrow Y$,
necessarily of degree $0$, are {\it isomorphisms\/} in this category.
Henceforth when we refer to isomorphisms in the multi setting
this kind of isomorphism will always be understood.

\proclaim{Proposition 1.7.3}
Let $X$ and $Y$ be bigraded $R$-modules, let $f \colon X \longrightarrow Y$
be a multimorphism, and let $Cf \colon CX \longrightarrow CY$ be the
corresponding morphism of graded $R$-modules. Then $f$ is an isomorphism if
and only if $Cf$ is an isomorphism of graded $R$-modules. \qed
\endproclaim

\proclaim{Lemma 1.8}
A multimorphism $f \colon X \to Y$ of degree $0$ of the kind
$f = f^0 + f^1 + f^2 + \cdots$is an isomorphism if and only if $f^0$
is an isomorphism.
\endproclaim
\demo{Proof}
It is obvious that the condition is necessary. To see that it is also
sufficient, suppose that $f^0$ is an isomorphism, and let $g^0$ be its inverse.
To extend $g^0$ to an inverse of $f$, all we have to do is to solve the equation
$$
\roman {Id} = f g = \sum f^i g^j
$$
for $g^1,g^2, \dots$
which amounts to solving the series
$$
0 = \sum_{i+j = k} f^i g^j, \quad k \geq 1,
$$
of equations for $g^1,g^2, \dots$. This series of equations
admits a unique solution $g^1,g^2, \dots$. \qed \enddemo

\noindent
{\smc Definition 1.9.1.}
Let $X$ and $Y$ be multicomplexes. A {\it morphism of multicomplexes\/}
written as $f \colon X \longrightarrow Y$, is a multimorphism
$$
f= f^0 + f^1 + f^2 + \cdots \colon X \longrightarrow Y
$$
of the underlying bigraded $R$-modules of degree zero having the property that
$$
df + (-1)^{|f|} fd = 0
$$
as multimorphisms of the underlying bigraded $R$-modules.

Thus in particular a morphism of multicomplexes preserves column filtrations.
\proclaim{Proposition 1.9.2}
Let $X$ and $Y$ be multicomplexes, and let $f \colon CX \longrightarrow CY$
be a morphism of chain complexes that preserves column filtrations.
Then the components
$$
f^k
\colon
X_{p,q}
\longrightarrow
X_{p-k,q+k},\quad k = 0,1,\cdots
$$
constitute a morphism of multicomplexes which, in turn, induces the original
morphism $f \colon CX \longrightarrow CY$ of filtered chain complexes
\endproclaim

\noindent
{\smc Definition 1.10.1}
Given two morphisms $f,g \colon X \longrightarrow Y$ of multicomplexes,
a {\it homotopy of morphisms of multicomplexes\/} or, more briefly, a
{\it multihomotopy\/}, written as
$$
h \colon f \simeq g \colon X \longrightarrow Y,
$$
is a multimorphism
$h= h^{-1} +h^0 + h^1 + h^2 + \cdots \colon X \longrightarrow Y$
of degree $1$ of the underlying bigraded $R$-modules satisfying
the identity
$$
dh + hd = g-f,
\tag1.10.2
$$
interpreted as one among multimorphisms of the underlying bigraded $R$-modules.
The two morphisms $f$ and $g$ of multicomplexes will
then be said to be {\it multihomotopic\/}.

Notice that a multihomotopy does not necessarily preserve column filtrations.
\proclaim{Proposition 1.10.3}
Let $X$ and $Y$ be multicomplexes, let $f,g \colon CX \longrightarrow CY$
be morphisms of chain complexes
that preserve
column filtrations,
and let
$h \colon CX \longrightarrow CY$
be a chain homotopy between
$f$ and $g$ that raises column filtration
at most by one.
Then the components
$$
h^k
\colon
X_{p,q}
\longrightarrow
X_{p-k,q+k+1},\quad k = -1,0,1,\cdots
$$
constitute a multihomotopy between $f$ and $g$, viewed as morphisms of
multicomplexes, and this multihomotopy induces the original chain homotopy
between $f$ and $g$, viewed as chain maps.
\endproclaim

\noindent
{\smc Definition 1.11.}
A {\it multiequivalence\/} is a morphism $f\colon X \longrightarrow Y$
of multicomplexes having an inverse with respect to the notion of multihomotopy;
in other words, $f$ is a multiequivalence, provided there are a morphism
$g\colon Y \longrightarrow X$ of multicomplexes and multihomotopies
$fg \simeq \roman{Id}$ and $gf \simeq \roman{Id}$.

\noindent {\smc Definition 1.12.1.} A {\it weak
multiequivalence\/} is a morphism $f\colon X \longrightarrow Y$ of
multicomplexes inducing an isomorphism
$$
f_* \colon \roman (E^r_{*,*}(X),d^r) \longrightarrow (E^r_{*,*}(Y),d^r)
$$
of column spectral sequences for $r \geq 2$.
\proclaim{Lemma 1.12.2}
Let $X$ and $Y$ be multicomplexes that are free as bigraded $R$-modules.
Then a weak multiequivalence $f\colon X \longrightarrow Y$
is a genuine multiequivalence.
\endproclaim
\demo{Proof} This is left to the reader. \qed \enddemo
Before we spell out the next observation we remind the
reader that
$\{F_p\}$
refers to the
column filtration reproduced in Definition 1.4.1 above.
\proclaim{Lemma 1.12.3}
A morphism $f\colon X \longrightarrow Y$ of multicomplexes is a weak
multiequivalence if and only if, for each $q \geq 0$, the restriction
$Cf|\colon F_q(CX) \longrightarrow F_q(CY)$ is a chain equivalence, that is,
if and only if $Cf\colon CX \longrightarrow CY$ is a filtered chain
equivalence with respect to the column filtration.
\endproclaim
\demo{Proof} This comes down to the standard identification of
$E^2(X)$ etc. with
\linebreak
$\Ho(F_*(CX), F_{*-1}(CX))$ etc.
Details are left to the reader. \qed \enddemo

\noindent
{\smc Definition 1.13.}
Let $X$ and $Y$ be multicomplexes. Then the {\it tensor product\/}
$X \otimes Y$ {\it in the category of multicomplexes\/}
is defined by
$$
\aligned
(X \otimes Y)_{p,q} &= \sum_{i+k=p,\ j+\ell=q} X_{i,j} \otimes Y_{k,\ell}
\\
d^r(x_{i,j} \otimes y_{k,l}) &=
d^r(x_{i,j}) \otimes y_{k,l} +
(-1)^{(i+j)}x_{i,j} \otimes d^r (y_{k,l}).
\endaligned
\tag 1.13.1
$$

Notice that when $C$ and $C'$ are chain complexes, the associated multicomplex
of their tensor product $C \otimes C'$ as chain complexes coincides with
the tensor product of the associated multicomplexes.

\noindent {\smc Definition 1.14.1.} The {\it horizontal
suspension\/}
 of a bigraded $R$-module $X$ is the bigraded $R$-module
$sX$ given by
$$
(sX)_{p,q} = X_{p-1,q};
\tag 1.14.2
$$
abusing notation somewhat, we write $s \colon X \to sX$ for  the corresponding
{\it (horizontal) suspension operator\/}, which is the identity when we
neglect bigrading and which, in the above language, is a multimorphism
of degree $\etad=1$ of the kind
$$
s =s^{\ell}\colon X_{*,*} \longrightarrow (sX)_{*+1,*},
$$
with $\ell = -1$, that is, $s$ has a single component.

\noindent
{\smc Definition 1.14.3}
The {\it suspension\/} of a multicomplex $X$ is the multicomplex
$sX$ which as a bigraded $R$-module is the horizontal suspension
and whose multidifferential is given by
$$
sd^j + d^j s = 0;
\tag 1.14.4
$$
here $s \colon X \to sX$ denotes the corresponding (horizontal) suspension
operator, and we do not distinguish in notation between the constituents of
the multidifferential on $X$ and $sX$.

Notice that when $C$ is a chain complex, the associated multicomplex
of its suspension $sC$ as a chain complex coincides with
the suspension of the associated multicomplex.

\noindent
{\smc Definition 1.15.}
Given two multicomplexes $X$ and $Y$, their {\it direct sum\/}
$X \oplus Y$ is the multicomplex given by
$$
(X \oplus Y)_{p,q} = X_{p,q}  \oplus Y_{p,q},
$$
with the obvious multidifferential induced by those on $X$ and $Y$.

\noindent {\smc Definition 1.16.} A {\it multialgebra\/} is a
bigraded algebra $A$ together with a multicomplex structure so
that the structure map $m \colon A \otimes A \to A$ is a morphism
of multicomplexes.

\noindent
{\smc Definition 1.17.}
Given a multialgebra $A$, a  {\it multi left\/} $A$-{\it module\/} $X$
is a multicomplex $X$ together with the structure
$m \colon A \otimes X \to X$
of a left bigraded $A$-module on $X$ that is a morphism of multicomplexes.
Multi right $A$-modules are defined accordingly.

Let $A$ and $B$ be  multialgebras, and let $f_1$ and $f_2$ be morphisms
$A \to B$ of multialgebras. Then $B$ admits an obvious structure
of a bigraded $A$-bimodule which we write as $(a,b) \mapsto a \cdot b$ where
$$
a \cdot b = (f_1(a))b,\quad b \cdot a = b (f_2(a)),
\quad a \in A,\, b \in B.
\tag1.17.1
$$
We shall refer to a multimorphism $h \colon A \to B$ of the underlying
bigraded $R$-modules as an $f_1$-$f_2$-{\it multiderivation\/},
provided it is a derivation with respect to the
bigraded $A$-bimodule structure ${\roman (1.17.1)}$, i.~e. if
$$
m (f_1 \otimes h + h \otimes f_2) =
h \,m,
\tag1.17.2
$$
where $m$ refers to the structure maps.

\noindent
{\smc Definition 1.18.}
A {\it homotopy\/} $f_1 \simeq f_2$  {\it of morphisms of multialgebras\/}
is a multihomotopy $h \colon A \to B$ (in the sense of {\rm (1.10)}) that is
also a $f_1$-$f_2$-multiderivation. More briefly we shall refer to such a
homotopy as a {\it multihomotopy\/} (in the context of morphisms of
multialgebras).

\noindent
{\smc 1.19.}
Given an ordinary differential graded algebra $A$, viewed as an ordinary chain
complex, the associated multicomplex (1.3) plainly inherits a multialgebra
structure which we refer the {\it associated multialgebra\/} structure.

\beginsection 2. Free multialgebras

\noindent {\smc Definition 2.1.} A multialgebra $A$ is  {\it
free\/} if its underlying bigraded algebra is (isomorphic to) the
tensor algebra $\roman T[V]$ on some free bigraded $R$-module $V$,
with the obvious bigrading, cf. (1.13).

A free multialgebra $A$ admits an obvious {\it augmentation map\/}
$\varepsilon \colon \roman T[V] \to R$, and we shall say it is
{\it connected\/} (as an augmented algebra) if $CV$ is
non-negative or if $CV$ is non-positive and zero in degree zero.

For convenience we recollect some properties of free connected multialgebras.
Henceforth $V$, $V'$, and $W$ denote free connected bigraded $R$-modules,
that are non-negative or non-positive and zero in degree zero,
and free bigraded algebras will always be assumed connected.

\noindent {\smc (2.2) Multiderivations and multidifferentials.}
Let $A=\roman T[V]$ be a free connected multialgebra, and let $M$
be a bigraded $A$-bimodule. As in Section 1 above, we refer to a
multimorphism
$$
d = d^0 + d^1 + d^2 + \cdots \colon
\roman T[V] \to \roman M
$$
of degree $-1$ that is also a derivation (with respect to the bigraded
$A$-bimodule structure) as a {\it multiderivation\/}. Each multiderivation
is plainly determined by its restriction
$$
\beta = d|\colon V \longrightarrow \roman M
$$
to $V$, and $\beta$ is a multimorphism. When $M$ itself is the bigraded
tensor algebra $\roman T[W]$ on some bigraded $R$-module $W$, the multimorphism
$\beta$ has {\it components\/}
$$
\beta_i \colon V \longrightarrow \roman W^{\otimes i}
$$
which are itself multimorphisms and, conversely, each sequence $\{\beta_i\}$
of multimorphisms $V \longrightarrow \roman W^{\otimes i}$ determines a
multiderivation $\roman T[V] \to \roman T[W]$.

In case $W=V$, for each multiderivation $d$ of degree $-1$, the
composite $dd$ is a multiderivation $d$ of degree $-2$, whence
$dd=0$ if and only if the restriction $dd|$ to $V$ vanishes. Hence
if $(\roman T[V],d)$ is a multialgebra, $(V,\beta_1)$ is a
multicomplex.

\noindent {\smc (2.3) Multimorphisms.}
Let $A=\roman T[V]$ be a
free multialgebra, and let $B$ be a bigraded algebra. Each
multimorphism $f \colon \roman T[V] \to  B$ of bigraded algebras
is determined by its restriction
$$
\alpha = f|\colon V \longrightarrow B.
$$
When $B$ is the bigraded tensor algebra $\roman T[W]$ on some bigraded
$R$-module $W$, the multimorphism $\alpha$ has {\it components\/}
$$
\alpha_i \colon V \longrightarrow \roman W^{\otimes i},
$$
and, conversely, each sequence $\{\alpha_i\}$ of multimorphisms
$V \longrightarrow \roman W^{\otimes i}$
determines a multimorphism $\roman T[V] \to \roman T[W]$ of bigraded algebras.

When $d$ and $d'$ endow $\roman T[V]$ and $\roman T[V']$, respectively, with
multialgebra structures and when $f \colon \roman T[V] \to\roman T[V']$
is a multimorphism of bigraded algebras,
$Df = d'f - fd$ is a multiderivation of degree $-1$
(with respect to the obvious bigraded
$\roman T[V]$-bimodule structure on $\roman T[V']$).
Hence $d'f = fd$ if and only if $d'f|V = fd|V$. In degree $1$ this condition
gives $\alpha_1 \beta_1 = \beta'_1 \alpha'_1$. Consequently if $f$ is a
morphism of multicomplexes, so is $\alpha_1$.

The proof of the following is straightforward and left to the
reader.
\proclaim{Lemma 2.3.1} A multimorphism $f \colon \roman
T[V] \to\roman T[V']$ of bigraded algebras is an isomorphism if
and only if its first component $\alpha_1\colon \roman V \to\roman
V'$ is an isomorphism. \qed
\endproclaim

\noindent {\smc (2.4) Multihomotopies.} Let $A$ and $B$ be
multialgebras, and let $f$ and $f'$ be morphisms $A \to B$ of
multialgebras. Recall from Section 1 that an
$f$-$f'$-multiderivation $h \colon A \to B$ of degree 1 is called
a {\it multihomotopy\/} $f \simeq f'$  {\it of morphisms of
multialgebras\/} provided $Dh (=dh + hd) = f' - f$. When  $A$,
viewed as a bigraded algebra, is a tensor algebra on some bigraded
$R$-module, this notion of multihomotopy can be conveniently
described in terms of a suitable cylinder construction, cf. e.~g.
\cite\baueslem\ and  \cite\huebkade\ (3.3)

Let $V$ be a multicomplex and let $A = \roman T[V]$. Then the {\it cylinder\/}
$A \times I$ is characterized as follows:
\roster
\item"--"
As a bigraded algebra, $A \times I$ is the tensor algebra
$\roman T[V' \oplus V'' \oplus sV]$
on the direct sum of two copies $V'$ and $V''$ of $V$ and the
horizontal suspension $sV$ of $V$ --
this is just the tensor algebra on the
corresponding cylinder $V \times I$;
we write $i' \colon A \to A \times I$
and
$i'' \colon A \to A \times I$
for the obvious injections of bigraded algebras
which identity $V$ with $V'$ and $V''$ respectively;
\item"--"
up to the obvious change in notation,
the multidifferential on $V'$ and $V''$
is the same as that in $V$;
\item"--"
to define the multidifferential on $sV$, let
$
S \colon
\roman T[V]
\longrightarrow
\roman A \times I
$
be the $i'$-$i''$-multiderivation determined by
$Sv=sv$,
so that, for $a,b \in \roman T[V]$
$$
S(ab) = (S(a))b'' + (-1)^{|a|}a'S(b),
\tag2.4.1
$$
and define the multidifferential $d$ on $sV \subseteq A \times I$
by
$$
\aligned
d^0(sv) &= -Sd^0v
\\
d^1(sv) &= v'' - v' -Sd^1v
\\
d^j(sv) &= -Sd^jv, \quad j \geq 2.
\endaligned
\tag2.4.2
$$
\endroster
We note that (2.4.2) implies that $S$ is a multihomotopy $ S
\colon i' \simeq i'' $ of morphisms of multialgebras. Moreover,
$i'$ and $i''$ are multiequivalences, and it is manifest that the
module of indecomposables $Q (A \times I)$ is just the
corresponding cylinder on the indecomposables $QA$.
\proclaim{Proposition 2.4.3} Let $V$ be a multicomplex, let $A =
\roman T[V]$, let $B$ be a multialgebra, and let $f'$ and $f''$ be
morphisms $A \to B$ of multialgebras. Then the formulas
$$
H i' = f',\quad H i'' = f'',
\quad
HS = h
$$
determine a natural bijection between multihomotopies
$
h \colon f' \simeq f''
$
of morphisms of multialgebras and morphisms
$
H \colon A \times I \longrightarrow B
$
of multialgebras with the property that
$$
H i' = f',\quad H i'' = f''.
$$
\endproclaim
\demo{Proof} This is straightforward and left to the reader. \qed
\enddemo

\noindent {\smc Remark 2.4.4.} It is not hard to deduce from
(2.4.3) that the above notion of multihomotopy of morphisms
between multialgebras $A$ and $B$ is an equivalence relation,
provided the $A$ underlying bigraded algebra is a tensor algebra;
cf. e.~g. \cite\baueslem\ for the more conventional case of chain
algebras. However, for arbitrary multialgebras $A$, this need {\it
not\/} be the case.

Let
$
H \colon \roman T[V] \times I \longrightarrow B
$
be a morphism of multialgebras. Then
$$
H(v'') = H(v') + Hdsv + HS(dv). \tag2.4.5
$$
Inspection shows that $H$ is determined by its values on $V'$ and $sV$.
In fact, write the multidifferential on $\roman T[V]$ in the form
$d = d^{(0)}  + \partial$ so that $d^{(0)}$ is the multidifferential that
comes from the multidifferential on $V$ and so that the
\lq\lq multi\rq\rq operator $\partial$ lowers augmentation filtration.
We note that this filtration has nothing to do with the corresponding row
or column filtrations; however, $\partial$ corresponds to a perturbation
of the differential induced by $d^{(0)}$ on the total complex
$C\roman T[V]=\roman T[CV]$ with respect to the augmentation filtration. Then
$$
Sdv = S d^0v + S\partial v = (d^0 v)' + (d^0 v)'' + \text {terms
of lower filtration}
\tag2.4.6
$$
whence by induction on degree we see that $H$ is determined by its
values on $V'$ and $sV$. This proves the following:
\proclaim{Lemma 2.4.7} Let $f \colon \roman T[V] \to B$ be a
morphism of multialgebras, and let $ \gamma\colon V \to \roman B $
be a multimorphism of degree 1 of the underlying bigraded modules.
Then there is a  morphism $ f' \colon \roman T[V] \to B $ of
multialgebras and a multihomotopy $h_{\gamma} \colon f \simeq f'$
of morphisms of multialgebras so that the restriction of
$h_{\gamma}$ to $V$ coincides with $\gamma$, and $f'$ and
$h_{\gamma}$ are uniquely determined by the given data. \qed
\endproclaim
\proclaim{Corollary 2.4.8} Let $f \colon \roman (T[V],d) \to
\roman (T[W],d)$ be a morphism between free multialgebras, let $
\alpha_1 \colon (V,\beta_1) \to (W,\tilde \beta_1) $ be its first
component, and let $ \alpha'_1 \colon (V,\beta_1) \to (W,\tilde
\beta_1) $ be a morphism of multicomplexes that is multihomotopic
to $\alpha_1$. Then there is a morphism
$$
g \colon \roman (T[V],d) \to \roman (T[W],d)
$$
of multialgebras with first component $\alpha'_1$ and multihomotopic to $f$.
\endproclaim

\proclaim{Theorem 2.5} {\rm (Multi Version of the Adams-Hilton
Theorem)} Let $A$ and $A'$ be multialgebras, not necessarily free
as bigraded modules over the ground ring, let
$$
\CD
@.
A
\\
@.
@VVgV
\\
(\roman T[V],d)
@>>{f'}>
A'
\endCD
\tag2.5.1
$$
be a diagram in the category of multialgebras, and suppose that
$g$ is a weak multiequivalence. Then there is a morphism
$
f
\colon
(\roman T[V],d)
@>>>
A$
of multialgebras so that $gf$ is homotopic to $f'$ as morphisms of multialgebras
and, furthermore, the multihomotopy class of $f$ is uniquely determined by this
condition.
\endproclaim
Another way to spell this out is to say that, for each weak multiequivalence
$g$, the induced morphism
$$
g^{\sharp} \colon [(\roman T[V],d), A] \longrightarrow [(\roman
T[V],d), A'] \tag2.5.2
$$
on the sets of homotopy classes of morphisms of multialgebras is a bijection.

A proof of the corresponding classical Theorem of Adams and Hilton
may be found in  \cite\adamhilt\ (3.1); see also \cite\baueslem\
(1.4). \demo{Proof of Theorem 2.5} For intelligibility, we
reproduce first the argument for the classical Theorem of Adams
and Hilton:

In degree $0$, the restriction $f| \colon \roman (T[V])_0 = R\to A_0$
is taken to be the obvious morphism that sends $1 \in R$ to $1 \in A$, and
the morphism $f \colon \roman T[V] \to A$ and  homotopy
$h \colon \roman T[V] \to A$ are then constructed by induction on the degree
of the generating module $V$. More precisely, appropriate  morphisms
$f_j\colon V_j \to A_j$
and
$h_j\colon V_j \to A'_{j+1}$
are constructed by induction in such a way that, for $j \geq 1$,
$$
dh_j + h_{j-1} d = g_j f_j - f'_j . \tag2.5.3
$$
Here are the details; we write $Z_k(-)$ etc. for cycles in degree $k$:
Let
$$
\zeta_1 \colon V_1 \to Z_1(A)
$$
be a
morphism so that
$$
g_1 f_1 - f'_1 \colon V_1 \longrightarrow Z_1(A')
$$
goes into the boundaries.
The existence of such a morphism $\zeta_1$ is guaranteed by the
hypothesis that $g$ is a weak equivalence.
Let
$$
f_1=\zeta_1 \colon V_1 \to Z_1(A).
$$
Since $V_1$ is free, there is a morphism
$$
h_1 \colon V_1 \to A'_2
$$
so that
$$
dh_1 = g_1 f_1 - f'_1 \colon V_1 \longrightarrow Z_1(A').
$$

Next, let $n \geq 1$, and suppose by induction that the components
$f_1,\dots,f_n$ and $h_1,\dots,h_n$ have already been constructed
in such a way that (2.5.3) holds for $1 \leq j \leq n$. Then the
composite
$$
g_n f_n d \colon V_{n+1} \longrightarrow A'_n
$$
goes into the $n$-cycles $Z_n(A')$ of $A'$, in fact, in view of
(2.5.3), we have
$$
g_n f_n d = f'_nd + dh_n d = d(f'_{n+1} + h_nd),
$$
whence
$g_n f_n d$
goes into the $n$-boundaries
of $A'$.
But $g$ is a weak equivalence,
and hence
$f_n d$
goes into the $n$-boundaries
of $A$, that is, there is a morphism
$$
\tilde f_{n+1}\colon V_{n+1} \longrightarrow A_{n+1}
$$
so that
$$
d\tilde f_{n+1} = f_n d.
$$
Moreover,
$$
\align
d(g_{n+1} \tilde f_{n+1} - f'_{n+1} -h_n d )&=
g_n f_n d - f'_n d - d h_n d
\\
 &=  d h_n d- d h_n d
\\
&= 0 .
\endalign
$$
Hence
$$
g_{n+1} \tilde f_{n+1} - f'_{n+1} -h_n d
\colon
V_{n+1} \longrightarrow
A'_{n+1}
$$
goes into the $(n+1)$-cycles
$Z_{n+1}(A')$ of $A'$.
Since $g$ is a weak equivalence,
there is a morphism
$$
\zeta_{n+1} \colon V_{n+1} \to Z_{n+1}(A)
$$
so that, with
$
f_{n+1} = \tilde f_{n+1} + \zeta_{n+1},
$
the morphism
$$
g_{n+1}  f_{n+1} - f'_{n+1} -h_n d
\colon
V_{n+1} \longrightarrow
A'_{n+1}
$$
goes into the $(n+1)$-boundaries
$Z_{n+1}(A')$ of $A'$.
Since $V_{n+1}$ is free, there is thus a morphism
$$
h_{n+1} \colon V_{n+1} \to A'_{n+2}
$$
so that
$$
dh_{n+1} + h_n d = g_{n+1} f_{n+1} - f'_{n+1} .
$$
This completes the inductive step.

Proceeding thus, as $n$ tends to infinity,
we obtain the desired morphism $f$ and homotopy $h$.

We now explain the necessary modifications for a complete argument
for Theorem 2.5, the multi version of the Adams-Hilton Theorem. We
shall show that, in the situation of Theorem 2.5, the morphism
$f\colon C\roman T[V] \to CA$ of differential graded algebras can
be constructed compatibly with the column filtrations and,
furthermore, that the homotopy $h$ can be constructed so that it
raises column filtration by 1. By Proposition 1.9.2, the morphism
$f$ then determines a corresponding morphism of multialgebras and,
by Proposition 1.10.3, the homotopy then determines a
corresponding multihomotopy. Here are the details:

As before, we denote the column filtrations by $\{F_q\}$.
For each $q \geq 0$, we then have
a diagram
$$
\CD
@.
F_q(A)
\\
@.
@VV{g|F_q(A)}V
\\
(\roman T[F_q(V)],d)
@>>{f'}>
F_q(A')
\endCD
$$
in the category of differential algebras; furthermore, cf. (1.12.3),
since $g$ is a weak multiequivalence, for each $q$, the restriction
$g|F_q(A)$ is a weak equivalence. By the corresponding classical
Adams-Hilton Theorem \cite\adamhilt, for each $q$, there is a morphism
$
f^{(q)}
\colon
(\roman T[F_q(V)],d)
@>>>
F_q(A)$
of differential graded algebras so that $(g|)f^{(q)}$ is homotopic to $f'|$
as morphisms of differential graded algebras and, furthermore, the homotopy
class of $f^{(q)}$ is uniquely determined by this condition. It remains to
show that the morphisms $f^{(q)}$ and the corresponding homotopies can be
constructed compatibly with the column filtrations. This is seen by a
slightly more complicated induction than the one that came into play above.
Here are the details for the inductive step.

Let $q \geq 0$, and suppose that the morphism
$$
f^{(q)} \colon F_q(\roman T[V])  = \roman T[F_q(V)]
\longrightarrow F_q(A)
$$
and chain homotopy
$$
h^{(q)} \colon F_q(\roman T[V]) = \roman T[F_q(V)]
\longrightarrow F_q(A)
$$
have already been constructed. Furthermore, let $n \geq 1$, and
suppose by induction that the components
$f^{(q+1)}_1,\dots,f^{(q+1)}_n$ and
$h^{(q+1)}_1,\dots,h^{(q+1)}_n$ have already been constructed in
such a way that the appropriate replacement for (2.5.3) holds for
$1 \leq j \leq n$, that is, that
$$
dh^{(q+1)}_j + h^{(q+1)}_{j-1} d = g^{(q+1)}_j f^{(q+1)}_j -
(f')^{(q+1)}_j \tag2.5.4
$$
for $1 \leq j \leq n$.
Then the composite
$$
g^{(q+1)}_n f^{(q+1)}_n d
\colon (F_{q+1}(V))_{n+1} \longrightarrow (F_{q+1}(A'))_n
$$
goes into the $n$-cycles $Z_n(F_{q+1}(A'))$ of $F_{q+1}(A')$, in
fact, in view of (2.5.4), we have
$$
g^{(q+1)}_n f^{(q+1)}_n d
= (f')^{(q+1)}_nd + dh^{(q+1)}_n d
= d((f')^{(q+1)}_{n+1} + h^{(q+1)}_nd),
$$
whence
$g^{(q+1)}_n f^{(q+1)}_n d$
goes into the $n$-boundaries
of $F_{q+1}(A')$.
But $g$ is a weak equivalence
which,
by virtue of (1.12.3),
is compatible with the filtrations
and hence
$f^{(q+1)}_n d$
goes into the $n$-boundaries
of $F_{q+1}(A)$, that is to say,
the morphism
$$
f^{(q)}_{n+1}\colon (F_{q}(V))_{n+1}
\longrightarrow F_{q}(A_{n+1})
$$
admits an
extension
$$
\tilde f^{(q+1)}_{n+1}\colon (F_{q+1}(V))_{n+1}
\longrightarrow F_{q+1}(A_{n+1})
$$
so that
$$
d\tilde f^{(q+1)}_{n+1} = f^{(q+1)}_n d.
$$
Moreover,
$$
\align
d(g^{(q+1)}_{n+1} \tilde f^{(q+1)}_{n+1}
- (f')^{(q+1)}_{n+1} -h^{(q+1)}_n d )&=
g^{(q+1)}_n f^{(q+1)}_n d
- (f')^{(q+1)}_n d - d h^{(q+1)}_n d
\\
 &=  d h^{(q+1)}_n d- d h^{(q+1)}_n d
\\
&= 0 .
\endalign
$$
Hence
$$
g^{(q+1)}_{n+1} \tilde f^{(q+1)}_{n+1} -
(f')^{(q+1)}_{n+1} -h^{(q+1)}_n d
\colon
(F_{q+1}(V))_{n+1} \longrightarrow
(F_{q+1}(A'))_{n+1}
$$
goes into the $(n+1)$-cycles $Z_{n+1}(F_{q+1}(A'))$ of
$F_{q+1}(A')$. Since $g$ is a filtered weak equivalence, cf. what
was said above, there is a morphism
$$
\zeta^{(q+1)}_{n+1} \colon (F_{q+1}(V))_{n+1} \to Z_{n+1}(F_{q+1}(A))
$$
so that, with
$
f^{(q+1)}_{n+1} = \tilde f^{(q+1)}_{n+1} + \zeta^{(q+1)}_{n+1},
$
the morphism
$$
g^{(q+1)}_{n+1}  f^{(q+1)}_{n+1}
- (f')^{(q+1)}_{n+1} -h^{(q+1)}_n d
\colon
(F_{q+1}(V))_{n+1} \longrightarrow
(F_{q+1}(A'))_{n+1}
$$
goes into the $(n+1)$-boundaries
$Z_{n+1}(F_{q+1}(A'))$ of $F_{q+1}(A')$.
Since $F_{q+1}(V_{n+1})$ is free,
the morphism
$$
h^{(q)}_{n+1} \colon F_{q}(V_{n+1}) \to (F_{q}(A'))_{n+2}
$$
admits an
extension
$$
h^{(q+1)}_{n+1} \colon F_{q+1}(V_{n+1}) \to (F_{q+1}(A'))_{n+2}
$$
so that
$$
dh^{(q+1)}_{n+1} + h^{(q+1)}_n d
= g^{(q+1)}_{n+1} f^{(q+1)}_{n+1} - (f')^{(q+1)}_{n+1} .
$$
This completes the inductive step.

Proceeding thus, as $n$ tends to infinity,
we obtain the desired extensions $f^{(q+1)}$ and $h^{(q+1)}$.
Likewise,
as $q$ tends to infinity,
we obtain the desired filtered
morphism $f$ and homotopy $h$.
This completes the proof.\qed
\enddemo
\beginsection 3. Minimal free multialgebras

Let $R$ be a local ring, with maximal ideal $\fra m \subseteq R$
and residue field $k$
and let $M$ be an $R$-module.
Recall that
a free resolution
$$
0
@<<<
M
@<{\varepsilon}<<
F_0
@<{\delta_1}<<
F_1
@<{\delta_2}<<
\cdots
\tag3.1
$$
of $M$ in the category of $R$-modules is called {\it minimal\/} if
$\delta_j(F_j) \subseteq \fra m\,F_{j-1}$ for $j \geq 1$, cf.
\cite\eilenber, \cite\serre. In particular, a minimal resolution
exists and is unique up to a (non-canonical) isomorphism of chain
complexes \cite\eilenber. We also recall the following

\noindent
{\smc Definition 3.2.}
A chain complex $(V,d)$ over (a local ring) $R$ that is free as a graded
module over $R$ (as always) is {\it minimal\/} provided
$d(V) \subseteq \fra m\,V$.

For us the key concepts will be those given in (3.3) and (3.5) below.

\noindent
{\smc Definition 3.3.}
A {\it minimal free multicomplex \/} is a multicomplex $(X,d)$ that is free
as a bigraded module over $R$ (as always) and has the properties that
$d^0 = 0$ and that $(X,d^1)$ is a minimal chain complex.

We note that a minimal free chain complex $(V,d)$ over $R$ is a minimal free
multicomplex with respect to the obvious multicomplex structure on $(V,d)$.
\proclaim{Proposition 3.4}
Let $(V,d)$ and $(V',d')$ be minimal free multicomplexes of finite type.
Then a weak multiequivalence
$$
f \colon
(V,d)
\longrightarrow
(V',d')
$$
is an isomorphism of multicomplexes.
\endproclaim
\demo{Proof} Since $(V,d)$ and $(V',d')$ are assumed to be minimal
(free) multicomplexes, the vertical differentials $d^0$ and
$(d')^0$ are zero, the operations $d^1 $ and $(d')^1$ are
horizontal differentials, and the \lq\lq component\rq\rq\ $f^0$ of
$f$ is a chain map
$$
f^0 \colon
(V,d^1)
\longrightarrow
(V',(d')^1).
$$
Furthermore, $f^0$ induces an isomorphism on homology. However,
$
(V,d^1)
$ and
$(V',(d')^1)$
are minimal free chain complexes
of finite type
in the sense of (3.2) whence
$f^0$ is an isomorphism of chain complexes
and in particular admits an inverse $g^0$;
cf. \cite\huebkade\ (5.3).
By (1.8), the morphism $g^0$ extends to an inverse $g$ of $f$
in the category of multicomplexes.
\qed \enddemo

\noindent
{\smc Definition 3.5.}
A {\it minimal free multialgebra\/} over $R$ is a free multialgebra
$(\roman T[V],d)$ such that the generating multicomplex $(V,\beta_1)$
is a minimal free multicomplex.

\proclaim{Proposition 3.6}
Let $(\roman T [V],d)$ and $(\roman T [V'],d')$ be minimal free multi
$R$-algebras, assume that $V$ and $V'$ are of finite type, and let
$
f \colon
(\roman T [V],d)
\longrightarrow
(\roman T [V'],d')
$
be a multiequivalence of multi $R$-algebras. Then $f$ is homotopic to an
isomorphism
$
g \colon
(\roman T [V],d)
\longrightarrow
(\roman T [V'],d')
$
of multi $R$-algebras by a homotopy of morphisms of multi $R$-algebras.
\endproclaim
\demo{Proof}
Let
$$
f= f^0 + f^1 + f^2 + \cdots
\colon \roman T [V] \longrightarrow \roman T [V']
$$
be the given morphism of multialgebras,
let
$$
\alpha = f| \colon V \longrightarrow \roman T [V']
$$
be the restriction of
$f$ to $V$, and write
$$
\alpha_i \colon V \longrightarrow (V')^{\otimes i}
$$
for its components, in the category of bigraded $R$-modules.
For each $i \geq 1$, $\alpha_i$ is then itself a multimorphism
$$
\alpha_i = \alpha_i^0 + \alpha_i^1 + \alpha_i^2 + \cdots
\colon V \longrightarrow (V')^{\otimes i},
$$
with the notion of tensor product (1.13) understood.
Furthermore,
$f^0$ is a morphism
$$
f^0 \colon \roman (T [V],d^1) \longrightarrow \roman (T [V'],(d')^1)
$$
of differential graded algebras,
and
its
first component
$$
\alpha_1^0 \colon
(V,\beta_1^1)
\longrightarrow
(V',(\beta')_1^1)
$$
is a morphism of chain complexes. However, since $f$ is a multiequivalence
of multi $R$-algebras, $f^0$ is a chain equivalence; in view of
\cite\huebkade\  (3.2.2), this implies that the first component $\alpha_1^0$
is a chain equivalence, in fact, by virtue of \cite\huebkade\ (5.3),
$\alpha_1^0$ is an isomorphism of chain complexes. In view of (1.8),
the morphism
$$
\alpha_1 \colon
(V,\beta_1)
\longrightarrow
(V',\beta'_1)
$$
is an isomorphism of multicomplexes. By Corollary 2.4.8,
$\alpha_1$ can be extended to a morphism
$$
g \colon
(\roman T [V],d)
\longrightarrow
(\roman T [V'],d')
$$
of multialgebras in such a way that (i) its first component
coincides with $\alpha_1$ and (ii) the morphisms $f$ and $g$ are
homotopic as morphisms of multialgebras. By virtue of Lemma 2.3.1,
the morphism $g$ is an isomorphism since so is $\alpha_1$. \qed
\enddemo

\noindent
{\smc Definition 3.7.}
Let $C$ be a chain complex over $R$ that is free as a graded $R$-module
(as always). Then a {\it minimal free multimodel for\/} $C$ is a minimal
free multicomplex $U$ over $R$ together with a multiequivalence
$\alpha \colon U \longrightarrow C$;
here $C$ is viewed as a multicomplex in the obvious way.

\noindent {\smc Definition 3.8.} Let $A$ be an augmented
differential graded algebra over $R$ that is free as a graded
$R$-module (as always). Then a {\it minimal free multimodel for\/}
$A$ is a minimal free multialgebra $(\roman T[V],d)$ over $R$
together with a morphism $ g \colon (\roman T[V],d)
\longrightarrow A $ of multi $R$-algebras that is also a
multiequivalence; here $A$ is identified with its associated
multialgebra (cf. 1.19).

\proclaim{Theorem 3.9} Let $A$ be a differential graded
$R$-algebra that is free as a module over $R$. If $(\roman
T[V],d,g)$ and $(\roman T[V'],d',g')$ are minimal free multimodels
for $A$, and if $V$ and $V'$ are of finite type, then $(\roman
T[V],d)$ and $(\roman T[V'],d')$ are isomorphic multialgebras.
\endproclaim
\demo{Proof} By Theorem 2.5, there is a multiequivalence $ \phi'
\colon (\roman T[V],d) \longrightarrow (\roman T[V'],d'). $ In
view of Proposition 3.6, this multiequivalence is multihomotopic
to an isomorphism $ \phi \colon (\roman T[V],d) \longrightarrow
(\roman T[V'],d') $ of multi $R$-algebras. \qed
\enddemo
\proclaim{Theorem 3.10} Let $A$ be a connected augmented
differential graded $R$-algebra whose underlying graded $R$-module
is free. Then $A$ has a minimal free multimodel. When the homology
of $A$ is of finite type, a minimal free multimodel is unique up
to isomorphism of multialgebras.
\endproclaim
\demo{Proof}
For $k \geq 1$, pick a minimal resolution of $\Ho_k(J\roman B A)$ and
assemble these resolutions to a minimal chain complex $(F \Ho,\delta)$
of the kind (4.2.b) in \cite\huebkade. The argument for the proof
of \cite\huebkade\  (4.4) yields a multidifferential
$\vartheta = \{\delta, \vartheta^2, \vartheta^3, \dots\}$ on
$F \Ho$ and a multiequivalence between
$(F \Ho,\vartheta)$ and $J\roman B A$, and
hence $(F \Ho,\vartheta)$ is a minimal free multimodel for
$J\roman B A$;
the kind of reasoning used at this stage may be found in III.1 of
\cite\berikash\ (for the special case where the base space is a point).
The construction in \cite\huebkade\  (4.7) then yields a free multimodel
$$
(\roman T[s^{-1}F {\Ho}(J\roman BA)],d')
\to (\roman T[s^{-1}(J\roman BA)],d_{\Omega})
\to
A
$$
for $A$ which is minimal by construction.

The uniqueness of the minimal free multimodel follows from Theorem 3.9.\qed
\enddemo

\noindent {\smc Remark 3.11.} In the situation of (the proof of)
Theorem 3.10, the differential $d'$ endows $M=s^{-1}F\Ho(JBA)$
with the structure of a {\it s(trongly) h(omotopy)
a(ssociative)\/} coalgebra; such a structure is dual to that of an
$A(\infty)$-algebra introduced in \cite\stasheff\ and christened
{\it s(trongly) h(omotopy) a(ssociative)\/} algebra in
\cite\stashalp. In the special case where $C$ is the  coalgebra of
normalized chains on a simply connected space $X$ and $A=\Omega
C$, the cobar construction on $C$, Theorem 3.10 above yields a
minimal model for the chain algebra of the loop space on $X$. For
the special case where the ground ring is (local as above and) a
principal ideal domain, a minimal model was obtained in
\cite\huebkade. In the  even more special case where the ground
ring is that of the reals, such a result may be found in
\cite\chen; in fact, over a field our construction boils down to
that of {\smcc Chen}. Related models, not necessarily over a local
ring and not necessarily minimal, have been developed in
\cite\huebkade. A special case of the kind of models in
\cite\huebkade\ may be found in \cite\gugentwo. For the dual
situation, i. e. where, instead of $C$, a differential graded
algebra $A$ is considered, related models are given in
\cite\gugensta \ (not necessarily minimal ones) and in
\cite\kadeishv. More comments may be found in Section 2 of
\cite\huebkade\ (see the discussion before (2.3) in the quoted
reference).

\bigskip

\widestnumber\key{999}

\centerline{\smc References}

\medskip\noindent

\ref \no  \adamhilt
\by J. F. Adams and P. J. Hilton
\paper On the chain algebra of a loop space
\jour  Comm. Math. Helv.
\vol 20
\yr 1955
\pages  305--330
\endref

\ref \no  \avrahalp \by L. L. Avramov and S. Halperin \paper
Through the looking glass \paperinfo in: Proceedings of a
conference held a Stockholm, 1984 \jour Lecture Notes in
Mathematics, vol 1183 \publ Springer \publaddr
Berlin--Heidelberg--New York \yr 1986 \pages 1--27
\endref
\ref \no  \baueslem
\by H. J. Baues and J. M. Lemaire
\paper Minimal models in homotopy theory
\jour  Math. Ann.
\vol 225
\yr 1977
\pages  219--242
\endref
\ref \no  \berikash \by N. A. Berikashvili \paper Differentials of
a spectral sequence \paperinfo (Russian) \jour Proc. Tbil. Math.
Institut \vol 51 \yr 1976 \pages 1--105
\endref
\ref \no  \brown
\by R. Brown
\paper The twisted Eilenberg--Zilber theorem
\jour Celebrazioni Archimedee del Secolo XX, Simposio di topologia
\yr 1964
\endref
\ref \no  \cartan
\by H. Cartan
\paper Alg\`ebres d'Eilenberg-Mac Lane et homotopie
\paperinfo expos\'es 2--11
\jour S\'eminaire H. Cartan, 1954--55
\publ\'Ecole Normale Sup\'erieure
\publaddr Paris
\yr 1956
\endref
\ref \no  \carteile
\by H. Cartan and S. Eilenberg
\book Homological Algebra
\publ Princeton University Press
\publaddr Princeton
\yr 1956
\endref
\ref \no  \chen
\by K.T. Chen
\paper Iterated path integrals
\jour Bull. Amer. Math. Soc.
\vol 83
\yr 1977
\pages  831--879
\endref
\ref \no \degrmosu
\by P. Deligne, P. Griffiths, J. Morgan, and D. Sullivan
\paper Real homotopy theory of K\"ahler manifolds
\jour Inv. Math.
\vol 29
\yr 1975
\pages  245--274
\endref
\ref \no \doldone
\by A. Dold
\paper Zur Homotopietheorie der Kettenkomplexe
\jour Math. Ann.
\vol 140
\yr 1960
\pages  278--298
\endref
\ref \no \doldtwo
\by A. Dold
\book Halbexakte Homotopiefunktoren
\bookinfo Lecture Notes in Mathematics
 No. 12
\publ Springer
\publaddr Berlin--Heidelberg--New York
\yr 1966
\endref
\ref \no \eilenber
\by S. Eilenberg
\paper Homological dimension and syzygies
\jour Ann. of Math.
\vol 64
\yr 1956
\pages  328--336
\endref
\ref \no \eilenmac
\by S. Eilenberg and S. Mac Lane
\paper On the groups ${\Ho(\pi,n)}$.I.
\jour Ann. of Math.
\vol 58
\yr 1953
\pages  55--106
\moreref
\paper II. Methods of computation
\jour Ann. of Math.
\vol 60
\yr 1954
\pages  49--139
\endref
\ref \no \gugenhei
\by V.K.A.M. Gugenheim
\paper On the chain complex of a fibration
\jour Illinois J. of Mathematics
\vol 16
\yr 1972
\pages 398--414
\endref
\ref \no \gugentwo
\by V.K.A.M. Gugenheim
\paper On a perturbation theory for the homology of the loop space
\jour J. of Pure and Applied Algebra
\vol 25
\yr 1982
\pages 197--205
\endref

\ref \no \gugenlam
\by V.K.A.M. Gugenheim and L. Lambe
\paper Perturbation in differential homological algebra
\jour Illinois J. of Mathematics
\vol 33
\yr 1989
\pages 566--582
\endref

\ref \no \gulasta
\by V.K.A.M. Gugenheim, L. Lambe, and J.D. Stasheff
\paper Algebraic aspects of Chen's twisting cochains
\jour Illinois J. of Math.
\vol 34
\yr 1990
\pages 485--502
\endref

\ref \no \gulstatw
\by V.K.A.M. Gugenheim, L. Lambe, and J.D. Stasheff
\paper Perturbation theory in differential homological algebra. II.
\jour Illinois J. of Math.
\vol 35
\yr 1991
\pages 357--373
\endref

\ref \no \gugenmay
\by V.K.A.M. Gugenheim and J.P. May
\paper On the theory and applications of differential
torsion products
\jour Memoirs of the Amer. Math. Soc.
\vol 142
\yr 1974
\endref
\ref \no \gugenmil
\by V.K.A.M. Gugenheim and J. Milgram
\paper On successive approximations in
homological algebra
\jour Trans. Amer. Math. Soc.
\vol 150
\yr 1970
\pages  157--182
\endref
\ref \no \gugenmun
\by V.K.A.M. Gugenheim and H. J. Munkholm
\paper On the extended functoriality of Tor and Cotor
\jour J. of Pure and Applied Algebra
\vol 4
\yr 1974
\pages  9--29
\endref
\ref \no \gugensta
\by V.K.A.M. Gugenheim and J.D. Stasheff
\paper On perturbations and ${A_{\infty}}$-structures
\jour Bull. Soc. Math. Belgique
\paperinfo Festschrift in honor of G. Hirsch's 60'th
birthday, ed. L. Lemaire
\vol 38
\yr 1986
\pages 237--245
\endref

\ref \no \halperin
\by S. Halperin
\paper Lectures on minimal models
\jour  Memoires de la Soc. Math. de France
\vol 9/10
\yr 1983
\endref

\ref \no \halpstas \by S. Halperin and J.D. Stasheff \paper
Obstructions to homotopy equivalences \jour Advances in Math. \vol
32 \yr 1979 \pages 233--278
\endref
\ref \no \heller
\by A. Heller
\paper Homological resolutions of complexes with operators
\jour  Ann. of Math.
\vol 60
\yr 1954
\pages  283--303
\endref
\ref \no \hilton \by P. J. Hilton \book Homotopy theory and
duality \publ Gordon and Breach Science Publishers \publaddr New
York-London-Paris \yr 1965
\endref
\ref \no \habili
\by J. Huebschmann
\paper Perturbation theory and small models for the chains of
certain induced fibre spaces
\paperinfo Habilitationsschrift Universit\"at Heidelberg 1984
\finalinfo {\bf Zbl.} 576.55012
\endref

\ref \no \homotype
\by J. Huebschmann
\paper The homotopy type of $F\Psi^q$. The complex and symplectic
cases
\paperinfo
in: Applications of Algebraic $K$-Theory
to Algebraic Geometry and Number Theory, Part II,
Proc. of a conf. at Boulder, Colorado, June 12 -- 18, 1983
\jour Cont. Math.
\vol 55
\yr 1986
\pages 487--518
\endref

\ref \no \perturba
\by J. Huebschmann
\paper Perturbation theory and free resolutions for nilpotent
groups of class 2
\jour J. of Algebra
\yr 1989
\vol 126
\pages 348--399
\endref

\ref \no \cohomolo
\by J. Huebschmann
\paper Cohomology of nilpotent groups of class 2
\jour J. of Algebra
\yr 1989
\vol 126
\pages 400--450
\endref
\ref \no \modpcoho
\by J. Huebschmann
\paper The mod $p$ cohomology rings of metacyclic groups
\jour J. of Pure and Applied Algebra
\vol 60
\yr 1989
\pages 53--105
\endref

\ref \no \intecoho
\by J. Huebschmann
\paper Cohomology of metacyclic groups
\jour Trans. Amer. Math. Soc.
\vol 328
\yr 1991
\pages 1-72
\endref

\ref \no \abelian
\by J. Huebschmann
\paper Cohomology of finitely generated abelian groups
\jour L'Enseignement Math\'ematique
\vol 37
\pages 61--71
\yr 1991
\endref

\ref \no \berikas \by J. Huebschmann \paper Berikashvili's functor
$\Cal D$ and the deformation equation \paperinfo Festschrift in
honor of N. Berikashvili's 70th birthday \jour Proceedings of the
A. Razmadze Mathematical Institute \vol 119 \yr 1999 \pages 59--72
\finalinfo {\tt math.AT/9906032}
\endref

\ref \no \holomorp \by J. Huebschmann \paper On the cohomology of
the holomorph of a finite cyclic group, \jour J. of Algebra (to
appear) \finalinfo{\tt math.GR/0303015}
\endref

\ref \no \omni \by J. Huebschmann \paper Higher homotopies and
Maurer-Cartan algebras: quasi-Lie-Rine\-hart, Gerstenhaber, and
Batalin-Vilkovisky algebras \paperinfo to appear in: The Breadth
of Symplectic and Poisson Geometry, Festschrift in honor of A.
Weinstein's 60th birthday; J. Marsden and T. Ratiu, eds.; Progress
in Mathematics \publ Birkh\"auser Verlag \publaddr Boston $\cdot$
Basel $\cdot$ Berlin \yr 2004 \finalinfo{\tt math.DG/0311294}
\endref

\ref \no \huebkade
\by J. Huebschmann and T. Kadeishvili
\paper Small models for chain algebras
\jour Math. Z.
\vol 207
\yr 1991
\pages 245--280
\endref

\ref \no \huebstas \by J. Huebschmann and J. D. Stasheff \paper
Formal solution of the master equation via HPT and deformation
theory \jour Forum mathematicum 14 \yr 2002 \pages 847--868
\finalinfo{\tt math.AG/9906036}
\endref

\ref \no \husmosta
\by D. Husemoller, J. C. Moore, and J. D. Stasheff
\paper Differential homological algebra and homogeneous spaces
\jour J. of Pure and Applied Algebra
\vol 5
\yr 1974
\pages  113--185
\endref

\ref \no \kadeishv
\by T.V. Kadeishvili
\paper On the homology theory of fibre spaces
\jour Uspekhi Mat. Nauk.
\vol 35:3
\yr 1980
\pages 183--188
\moreref
\paperinfo translated in:
\jour Russian Math. Surveys
\vol 35:3
\yr 1980
\pages 231--238
\endref

\ref \no \kadeifou
\by T. Kadeishvili
\paper The predifferential of a twisted product
\jour Russian Math. Surveys
\vol 41
\yr 1986
\pages 135--147
\endref
\ref \no  \kadeifiv \by T. Kadeishvili \paper $A_{\infty}$-algebra
Structure in Cohomology and Rational Homotopy Type
 \paperinfo
(Russian) \jour Proc. Tbil. Math. Institut \vol 107 \yr 1993
\pages 1--94
\endref

\ref \no \lambstas
\by L. Lambe and J. D. Stasheff
\paper Applications of perturbation theory to iterated fibrations
\jour Manuscripta Math.
\vol 58
\yr 1987
\pages 363--376
\endref
\ref \no \liulevon
\by A. Liulevicius
\paper Multicomplexes and a general change of rings theorem
\paperinfo mi\-meo\-graphed notes, University of Chicago
\endref
\ref \no \liulevtw
\by A. Liulevicius
\paper A theorem in homological algebra and stable homotopy of projective
spaces
\jour Trans. Amer. Math. Soc.
\vol 109
\yr 1963
\pages  540--552
\endref

\ref \no \maclane
\by S. Mac Lane
\book Homology
\bookinfo Die Grundlehren der mathematischen Wissenschaften
 No. 114
\publ Springer
\publaddr Berlin--G\"ottingen--Heidelberg
\yr 1963
\endref
\ref \no \meyer
\by J. P. Meyer
\paper Acyclic models for multicomplexes
\jour Duke Math. J.
\vol 45
\yr 1978
\pages 76--85
\endref
\ref \no \mooreone
\by J. C. Moore
\paper Alg\`ebres d'Eilenberg-Mac Lane et homotopie
\paperinfo expos\'es 12 et 13
\jour S\'eminaire H. Cartan, 1954--55
\publ\'Ecole Normale Sup\'erieure
\publaddr Paris
\yr 1956
\endref
\ref \no \mooretwo
\by J. C. Moore
\paper Differential homological algebra
\jour Actes du Congr. Intern. des Math\'ematiciens
\yr 1970
\pages  335--339
\endref
\ref \no \moorethr
\by J. C. Moore
\paper Cartan's constructions
\paperinfo Colloque analyse et topologie, en l'honneur de Henri
Cartan
\jour Ast\'erisque
\vol 32--33
\yr 1976
\pages  173--221
\endref
\ref \no \moorsmit
\by J. C. Moore and L. Smith
\paper Hopf algebras and multiplicative fibrations
\jour Amer. J. of Math.
\vol 40
\yr 1968
\pages 752--780
\endref
\ref \no \munkholm
\by H. J. Munkholm
\paper The Eilenberg--Moore spectral sequence and strongly homotopy
multiplicative maps
\jour J. of Pure and Applied Algebra
\vol 9
\yr 1976
\pages  1--50
\endref
\ref \no \quillen
\by D. Quillen
\paper Rational homotopy theory
\jour Ann. of Math.
\vol 90
\yr 1969
\pages  205--295
\endref

\ref \no \sanebone
\by S. Saneblidze
\paper Homology classification of differential algebras
\jour Bulletin of the Academy of Sciences of the Georgian SSR
\vol 129
\yr 1988
\pages 241--243
\finalinfo  (Russian. Georgian summary)
\endref
\ref \no \sanebthr \by S. Saneblidze \paper Filtered model of a
fibration and rational obstruction theory \jour manuscripta math.
\vol 76 \yr 1992 \pages 111--136
\endref
\ref \no \sanebsev \by S. Saneblidze \paper The homotopy
classification of spaces by the fixed loop space homology
\paperinfo Festschrift in honor of N. Berikashvili's 70th birthday
\jour Proceedings of the A. Razmadze Mathematical Institute \vol
119 \yr 1999 \pages 155-164
\endref

\ref \no \serre
\by J. P. Serre
\book Alg\`ebre locale. Multiplicit\'es
\bookinfo Lecture Notes in Mathematics, No.~11
\publ Springer
\publaddr Berlin-Heidelberg-New York
\yr 1965
\endref
\ref \no \shih
\by W. Shih
\paper Homologie des espaces fibr\'es
\jour Pub. Math. Sci. IHES
\vol 13
\yr 1962
\endref
\ref \no \smirnov
\by V.A. Smirnov
\paper Homology of fibre spaces
\jour Russ. Math. Surveys
\vol 35
\yr 1980
\pages 294--298
\endref
\ref \no \stasheff
\by J.D. Stasheff
\paper Homotopy associativity of H-spaces.I
\jour Trans. Amer. Math. Soc.
\vol 108
\yr 1963
\pages 275--292
\moreref
\paper II
\jour Trans. Amer. Math. Soc.
\vol 108
\yr 1963
\pages 293--312
\endref
\ref \no \stashalp
\by J.D. Stasheff and S. Halperin
\paper Differential algebra in its own rite
\jour Proc. Adv. Study Alg. Top. August 10--23, 1970, Aarhus, Denmark
\pages 567--577
\endref
\ref \no \sullivan \by D. Sullivan \paper Differential forms and
the topology of manifolds \jour Proc. Conf. Manifolds Tokyo \yr
1973
\endref
\ref \no \sullitwo
\by D. Sullivan
\paper Infinitesimal Computations in Topology
\jour Pub. Math. I. H. E. S
\vol 47
\yr 1978
\pages  269--331
\endref
\ref \no \wall
\by C.T.C. Wall
\paper Resolutions for extensions of groups
\jour Proc. Camb. Phil. Soc.
\vol 57
\yr 1961
\pages 251--255
\endref

\enddocument